\documentclass[10pt,a4paper,twocolumn]{IEEEtran}

\IEEEoverridecommandlockouts                              
\usepackage{amsmath}
\usepackage{color}
\usepackage{graphicx}
\usepackage{latexsym}

\usepackage{epsfig} 

\title{{Minimal Inputs/Outputs for Subsystems in a Networked System}}

\author{Tong Zhou
\thanks{This work was supported in part by the NNSFC under Grant 61573209 and 51361135705.}
\thanks{T.Zhou is with the Department of Automation, Tsinghua University, Beijing, 100084,
CHINA. (Tel: 86-10-62797430; Fax: 86-10-62786911; e-mail:
tzhou@mail.tsinghua.edu.cn.)}}

\begin{document}
\renewcommand{\thefootnote}{\fnsymbol{footnote}}
\maketitle 
\renewcommand{\thepage}{10--\arabic{page}}
\setcounter{page}{1}

\begin{abstract}
Minimal input/output selection is investigated in this paper for {each subsystem of} a networked system. Some novel sufficient conditions are derived respectively for the controllability and observability of a networked system, as well as some necessary conditions. These conditions only depend separately on parameters of each subsystem and its in/out-degrees. It {is proven}  that in order to be able to construct a controllable/observable networked system, it is necessary and sufficient that each subsystem is controllable/observable. In addition, both sparse and dense subsystem connections are helpful in making the whole system controllable/observable. An explicit formula is given for the smallest number of inputs/outputs for each subsystem required to guarantee controllability/observability of the whole system.

\end{abstract}

\IEEEpeerreviewmaketitle

\section{Introduction}

With the increment of the dimension of a system, which is mainly due to technology developments in sensors, communications, etc., as well as more complicated and demanding tasks expected for a system, computation costs and numerical stability emerge as essential issues in system analysis and synthesis \cite{siljak78,sbkkmpr11,zhou15,zz16a}. It is now widely recognized that with the increment of its subsystem number, direct applications of results about a lumped system to a large scale or networked system may often result in an exponential increment of computation time and storage requirements \cite{hot11,lsb11,siljak78,sbkkmpr11,zhou15}. To make things worse, these direct applications are usually numerically unreliable. A well known example is to compute the eigenvalues and/or eigenvectors of a square matrix, which is often required in analyzing system properties and designing a controller. If this matrix has a large dimension and its condition number is also large, then numerical computation results of all the available algorithms are generally far from actual values \cite{hj91,zdg,zhou15}. To overcome these difficulties, various efforts have emerged recently for the analysis and synthesis of a networked system. Among which, an extensively studied problem is about its controllability/observability verifications, and construction of a controllable/observable networked system \cite{cd13,lsb11,Olshevsky14,pka16,ra13,scl15,trpj16,zhou15,zhou16a}.

Various results have now been obtained for this important theoretical issue on systems and control. For example, robustness of structural controllability, input addition, decentralized controllability, etc., have been investigated respectively in \cite{cd13,ra13}. In \cite{pzb14}, clustered networks are found easier to be controlled. It is declared in \cite{Olshevsky14} that finding the sparest input/output matrix such that a networked system is controllable/observable is NP-hard, and some algorithms are suggested in \cite{scl15,pka16} to approximately solve this minimal controllability/observability problem. A minimal actuator placement problem is also proven in \cite{trpj16} to be NP-hard, and a best approximation is suggested which has a polynomial computational complexity. Structural controllability and the cavity method are used in \cite{lsb11} to derive a set of driver nodes for assuring system controllability. In \cite{zhou15}, we have obtained a necessary and sufficient condition for an arbitrarily connected networked system to be controllable/observable, which depends separately on parameters of each subsystem. These results have been extended to various situations, such as the full column normal rank (FCNR) condition adopted in \cite{zhou15} is not satisfied, there are constraints on system inputs and states, etc. \cite{zz15,zhou16a,zz17}. {It has been discovered in \cite{sm68,yzdwl13,zhou17} that, when the state transition matrix (STM) of a networked system is given}, in order to guarantee its controllability/observability, the minimal number of its inputs/outputs is equal to the maximum geometric multiplicity of its STM.

In actual engineering, however, it is generally preferable to have inputs/outputs directly and separately affecting/measuring the states of each individual subsystem and/or their functions \cite{siljak78,scl15,pka16,wj01,zhou15}. Under this restriction, it is still not clear how many inputs/outputs are required for each subsystem to make the whole system controllable/observable. To emphasize this characteristic, the associated problem is called in this paper local input/output selections.

To settle this problem, we at first investigate relations among subsystem observability/controllabilty, subsystem out/in-degree and system observability/controllability. It has been made clear that in order to guarantee the observability/controllability of a networked system, each of its subsystems should be observable/controllable. A sufficient condition is derived for system observability which depends separately only on parameters of each subsystem and its out-degree. This condition reveals that both sparse and dense subsystem  connections are helpful to make the whole system observable/controllable. On the basis of these results, it is further proven that a necessary and sufficient condition for being able to build an observable/controllable networked system is that, each subsystem is observable/controllable. {It has also been proven} that in order to guarantee system controllability/observability, the number of inputs/outputs in each subsystem must at least be equal to that of the maximum geometric multiplicity of its STM.

The outline of this paper is as follows. At first, Section II gives a precise problem formulation and some preliminary results. Relations between controllability/observability of a networked dynamic system and subsystem out/in-degree are investigated in Section III. The minimal input/output problem is discussed in Section IV. Finally, some concluding remarks are given in Section V. An appendix is included to give proofs of some technical results. Some numerical examples are provided to illustrate the obtained theoretical conclusions.

The following notation and symbols are adopted. ${\cal R}^{m\times n}$ and ${\cal C}^{m\times n}$ are utilized respectively to represent the sets of $m\times n$ dimensional real and complex matrices. When $m$ and/or $n$ are equal to $1$, they are usually omitted. ${\rm\bf diag}\!\{X_{i}|_{i=1}^{L}\}$ denotes a block diagonal matrix with its $i$-th diagonal block being $X_{i}$, while ${\rm\bf col}\!\{X_{i}|_{i=1}^{L}\}$ the vector/matrix stacked by $X_{i}|_{i=1}^{L}$ with its $i$-th row block vector/matrix being $X_{i}$. $0_{m}$ and $0_{m\times n}$ represent respectively the $m$ dimensional zero column vector and the $m\times n$ dimensional zero matrix. The superscripts $T$ and $H$ stand respectively for the transpose and the conjugate transpose of a matrix/vector, while $||\cdot||_{2}$ the Euclidean norm of a vector.

\section{{Problem Formulation and Some Preliminaries}}

{Consider the networked system ${\rm\bf\Sigma}$ adopted in \cite{zhou15,zhou16a,zz15,zz17} which consists of $N$ linear time invariant (LTI) dynamic subsystems.} In this system, the dynamics of its $i$-th
subsystem ${\rm\bf\Sigma}_{i}$ is described by
\begin{equation}
\left[\!\!\begin{array}{c} x(t\!+\!1,i) \\ z(t,i) \\
y(t,i)
\end{array}\!\!\right] \!\!=\!\!\left[\!\!\begin{array}{ccc} A_{\rm\bf TT}(i) & A_{\rm\bf TS}(i) & B_{\rm\bf T}(i)  \\
A_{\rm\bf ST}(i) & A_{\rm\bf SS}(i) & B_{\rm\bf S}(i) \\
C_{\rm\bf T}(i) & C_{\rm\bf S}(i) & D(i)
\end{array}\!\!\right]\! \left[\!\!\begin{array}{c} x(t,i)
\\ v(t,i) \\ u(t,i)
\end{array}\!\!\right]  \label{eqn:1}
\end{equation}
and interactions among its subsystems are described by
\begin{equation}
v(t)=\Phi z(t) \label{eqn:2}
\end{equation}
Here, $z(t)={\rm\bf col}\!\left\{z(t,i)|_{i=1}^{N}\right\}$ and $v(t)={\rm\bf
col}\!\left\{v(t,i)|_{i=1}^{N}\right\}$. Moreover, $t$ and $i$ stand respectively for
the temporal variable and the index number of a subsystem, $x(t,i)$ represents the state vector of the
$i$-th subsystem ${\rm\bf\Sigma}_{i}$ at time $t$, $z(t,i)$ and
$v(t,i)$ respectively its outputs affecting other subsystems and
inputs denoting influences from other subsystems, $y(t,i)$ and $u(t,i)$
respectively its output and input vectors. Similar to the treatments adopted in \cite{zhou15,zhou16a}, in order to distinguish these vectors, $z(t,i)$ and $v(t,i)$ are called internal output/input vectors, while $y(t,i)$ and $u(t,i)$ external output/input vectors.

The following hypotheses are adopted throughout this paper.

\hspace*{-0.45cm}{\rm\bf Assumption 1:} the dimensions of the vectors $x(t,i)$, $v(t,i)$, $u(t,i)$, $z(t,i)$ and $y(t,i)$, are respectively $m_{{\rm\bf x}i}$, $m_{{\rm\bf v}i}$, $m_{{\rm\bf u}i}$, $m_{{\rm\bf z}i}$ and $m_{{\rm\bf y}i}$.  \hspace{\fill}$\Diamond$

\hspace*{-0.45cm}{\rm\bf Assumption 2:} the networked system ${\rm\bf \Sigma}$ is well-posed.  \hspace{\fill}$\Diamond$

\hspace*{-0.45cm}{\rm\bf Assumption 3:} the subsystem connection matrix (SCM) $\Phi$ is a constant matrix, and each of its rows  has only one nonzero element which is equal to one.  \hspace{\fill}$\Diamond$

The {first} assumption is adopted only for clarifying dimensions of the associated vectors, while the {second} one is necessary for a networked system to work properly, {which physically means that for an arbitrary external input series
${\rm\bf col}\!\left\{u(t,i)|_{i=1}^{N}\right\}_{t=0}^{\infty}$, the system states ${\rm\bf col}\!\left\{x(t,i)|_{i=1}^{N}\right\}_{t=0}^{\infty}$, as well as the external outputs ${\rm\bf col}\!\left\{y(t,i)|_{i=1}^{N}\right\}_{t=0}^{\infty}$, are uniquely determined} \cite{zdg,zhou15}. This assumption is equivalent to the requirement that the matrix $I-\Phi{\rm\bf diag}\!\left\{\!A_{\rm\bf SS}(i)|_{i=1}^{N}\!\right\}$ is invertible \cite{zhou15}. The {third} assumption appears very restrictive, but as argued in \cite{zhou15,zhou16a,zz15}, it actually does not introduce any constraints on the structure of the whole system. Briefly, when this assumption is not satisfied by an original system model, {it can be satisfied by a modified model with completely the same input-output relations,} through simply augmenting the associated subsystem internal input/output vectors $v(t,i)/z(t,i)$ with repeated elements, and modifying the associated matrices $A_{\rm\bf ST}(i)$, $A_{\rm\bf SS}(i)$ and $B_{\rm\bf S}(i)$. {Note that a large scale networked system usually has a sparse structure, which implies that this augmentation generally does not increase significantly the dimensions of the associated matrices.} In addition, under this assumption, each element of a subsystem's internal output vector is able to simultaneously affect more than one subsystems, and different elements of an internal output vector are able to affect different subsystems.

{In this paper, the following problem is investigated.}

\hspace*{-0.45cm}{{\rm\bf Problem:} For prescribed subsystem STMs $\left.A_{\rm\bf TT}(i)\right|_{i=1}^{N}$, find the minimal $m_{{\rm\bf z}i}$ and $m_{{\rm\bf y}i}$ ($m_{{\rm\bf v}i}$ and $m_{{\rm\bf u}i}$) , such that an observable (controllable) networked system ${\rm\bf\Sigma}$ can be constructed using only external outputs ${\rm\bf col}\!\left\{y(t,i)|_{i=1}^{N}\right\}$, $t=0, 1, 2,\cdots$ (external inputs ${\rm\bf col}\!\left\{u(t,i)|_{i=1}^{N}\right\}$, $t=0, 1, 2,\cdots$).  \hspace{\fill}$\Diamond$}

{A similar problem has been investigated in \cite{sm68,yzdwl13,zhou17} for a lumped system. The above problem, however, is different in the sense that it asks for the minimal number of outputs/inputs for {\it each} subsystem in constructing an observable/controllable networked system in the whole. This requirement reflects the fact  that subsystems of a networked system are usually far away from each other geometrically, which makes it expensive in engineering practices to have a signal that simultaneously and directly affects actuators of two or more different subsystems, or have a sensor to measure an output that is an explicit function of the states of several subsystems. In other words, it is more attractive in applications to restrict each input to {\it directly} affect states of {\it only} one subsystem, as well as to restrict each output to be a {\it direct} linear combination of the states {\it only} in one subsystem. To emphasize this characteristic, an input/output meeting these restrictions is called a {\it local} input/output, and the associated input/output selection problem is called a minimal {\it local} input/output problem.}

To investigate this problem, the following results are required which are widely known as the PBH test \cite{zdg,zhou15}.

\hspace*{-0.45cm}{\bf Lemma 1.} Consider a discrete LTI system with the following state space model
\begin{equation}
x(t+1)=Ax(t)+Bu(t),\hspace{0.5cm}
y(t)=Cx(t)+Du(t)\label{eqn:3}
\end{equation}
\begin{itemize}
\item This system is controllable, if and only if for every complex scalar $\lambda$,  the matrix $\left[\lambda I-A \;\; B\right]$ is of full row rank (FRR).
\item This system is observable, if and only if for every complex scalar $\lambda$,  the matrix ${\rm\bf col}\!\{\lambda I-A,\; C\}$ is of full column rank (FCR). \hspace{\fill}$\Diamond$
\end{itemize}

We sometimes also use an expression like that the matrix pair $(A,\:C)$ is observable, and that the matrix pair $(A,\:B)$ is controllable, when the associated system is.

The next lemma provides some characteristics of a plant transmission zero, which is closely related to the existence of a nonzero plant input vector sequence that makes its output vector constantly equal to zero \cite{zdg}.

\hspace*{-0.45cm}{\bf Lemma 2.} Let $G(\lambda)$ be a proper transfer function matrix having FCNR. Then, a complex number $\lambda_{0}$ is one of its transmission zeros,  if and only if there exists a nonzero complex vector $z_{0}$ satisfying $G(\lambda_{0})z_{0}=0$.  \hspace{\fill}$\Diamond$

To make mathematical derivations more concise, the following matrix symbols are adopted throughout this paper. $A_{\rm\bf
*\#}\!\!=\!\!{\rm\bf diag}\!\left\{\!A_{\rm\bf
*\#}(i)|_{i=1}^{N}\!\right\}$, $B_{\rm\bf *}\!\!=\!\!{\rm\bf
diag}\!\!\left\{\!B_{\rm\bf *}(i)|_{i=1}^{N}\!\right\}$,
$C_{\rm\bf *}\!\!=\!\!{\rm\bf diag}\!\!\left\{\!C_{\rm\bf
*}(i)|_{i=1}^{N}\!\right\}$ and
$D\!=\!{\rm\bf diag}\!\left\{D(i)|_{i=1}^{N}\right\}$, in which
${\rm\bf *,\#}={\rm\bf T}$ or ${\rm\bf S}$. Moreover, denote
${\rm\bf col}\!\left\{u(t,i)|_{i=1}^{N}\right\}$, ${\rm\bf
col}\!\left\{x(t,i)|_{i=1}^{N}\right\}$ and ${\rm\bf
col}\!\left\{y(t,i)|_{i=1}^{N}\right\}$ respectively by
$u(t)$, $x(t)$ and $y(t)$. Furthermore, define integers $M_{{\rm\bf \star}i}$ and $M_{\rm\bf \star}$ as $M_{\rm\bf
\star}={\sum_{k=1}^{N} m_{{\rm\bf \star}k}}$, $M_{{\rm\bf \star}0}=0$, $M_{{\rm\bf \star}i}={\sum_{k=1}^{i} m_{{\rm\bf \star}k}}$ with $1\leq
i\leq N$. Here, ${\rm\bf \star}={\rm\bf x}$, ${\rm\bf u}$, ${\rm\bf y}$, ${\rm\bf v}$ or ${\rm\bf z}$.

The following results have been established in \cite{zhou15} which are starting points of this paper.

\hspace*{-0.45cm}{\bf Lemma 3.} Define a matrix valued polynomial $M(\lambda)$ as
\begin{equation}
M(\lambda)=\left[\begin{array}{cc}
\lambda I_{M_{\rm\bf x}}-A_{\rm\bf TT} & -A_{\rm\bf TS} \\
-C_{\rm\bf T} & -C_{\rm\bf S} \\
-\Phi A_{\rm\bf ST} & I_{M_{\rm\bf v}}-\Phi A_{\rm\bf SS} \end{array}\right]     \label{eqn:4}
\end{equation}
The networked system ${\rm\bf\Sigma}$ is observable, if and only if at each complex scalar $\lambda$, $M(\lambda)$ is of FCR.  \hspace{\fill}$\Diamond$

The following results give the minimal number of outputs of a lumped system for observability assurance, which appears to be firstly observed in \cite{sm68} and re-observed in \cite{yzdwl13}. {Their correct proof, however, seems to be firstly given in \cite{zhou17}, in which the requirement that an output matrix must be real valued has been taken into account.}

\hspace*{-0.45cm}{{\bf Lemma 4.} Concerning the LTI system of Equation (\ref{eqn:3}), there exists a matrix $C$ such that this system is observable, if and only if the dimension of the output vector $y(t)$ is not smaller than the maximum geometric multiplicity of the STM $A$.               \hspace{\fill}$\Diamond$}

\vspace{-0.0cm}
\begin{table*}[!ht]
\vspace*{-0.00cm}
\caption{\vspace*{-0.00cm} Singular Values of the Observability Matrices for Each Subsystem}
\begin{center}
\vspace{-0.2cm}
\begin{tabular}{|p{2.0cm}|p{1.8cm}|p{1.8cm}|p{1.8cm}|p{1.8cm}|p{1.8cm}|}
\hline
 \multicolumn{3}{|c|}{External Outputs Only} & \multicolumn{3}{c|}{External + Internal Outputs} \\
\cline{1-6}
 Subsystem 1 & Subsystem 2 & Subsystem 3 & Subsystem 1 & Subsystem 2 & Subsystem 3 \\
 \hline
\hline
$8.4218\times 10^{-10}$ & $4.4831\times10^{-1}$ & $3.7206\times10^{-1}$ & $4.2694\times10^{-1}$ & $4.5542\times10^{-1}$ & $5.9149\times10^{-1}$ \\ \hline
$1.1603$               &	$1.0254$	       & $1.5623$              & $2.0086$              & $1.0361$              & $1.8373$ \\ \hline
$5.4545$               &	$6.0747$	       & $5.7118$              & $5.6423$              & $6.0810$              & $5.8053$ \\ \hline
\end{tabular}
\vspace{-0.5cm}
\end{center}
\end{table*}


\section{Out-degree, Controllability and Observability of a Networked System}

To investigate the minimal local input/output selection problem, we at first develop some new methods for verifying the controllability and observability of the networked system $\bf\rm\Sigma$. For this purpose, the following property of the SCM $\Phi$ is at first introduced. This property is firstly observed in \cite{zz16a} and plays important roles in the analysis of its stability and robust stability.

Let $m(i)$ stand for the number of subsystems that is directly affected by the $i$-th element of the vector $z(t)$, $i=1,\;2,\;\cdots,\;M_{\rm\bf z}$. Define matrices $\Theta(j)$, $j=1,2,\cdots,N$, and $\Theta$ respectively as $\Theta(j)={\rm\bf diag}\{\sqrt{m(i)}|_{i=M_{{\rm\bf z},j-1}+1}^{M_{{\rm\bf z},j}}\}$ and $\Theta={\rm\bf diag}\{\sqrt{m(i)}|_{i=1}^{M_{\rm\bf z}}\}$. It has been proven in \cite{zz16a} that
\begin{equation}
\Phi^{T}\Phi = \Theta^{2}
={\rm\bf diag}\left\{\left.\Theta^{2}(j) \right|_{j=1}^{N}\right\}
\label{eqn:5}
\end{equation}

Obviously from the definition of $m(i)$, we have that $\sum_{i=M_{{\rm\bf z},j-1}+1}^{M_{{\rm\bf z},j}}m(i)$ equals the out-degree of the $j$-th subsystem of the networked system $\rm\bf\Sigma$.

On the basis of this relation and Lemma 3, a necessary condition is obtained for the observability of System $\bf\rm\Sigma$. Its proof is given in the appendix.

\hspace*{-0.45cm}{\bf Lemma 5.} The networked system ${\rm\bf\Sigma}$ is observable, only if for each $i=1,2,\cdots,N$, the matrix pair $(A_{\rm\bf TT}(i),$ ${\rm\bf col}\!\{ C_{\rm\bf T}(i),\; A_{\rm\bf ST}(i)\})$ is observable.  \hspace{\fill}$\Diamond$

From the state space model of the subsystem ${\rm\bf\Sigma}_{i}$, it is clear that both the vector $y(t,i)$ and the vector $z(t,i)$ are its output vectors. In other words, when this subsystem is isolated from the influences of other subsystems, and its influences to other subsystems are also completely removed, then the observability of the subsystem ${\rm\bf\Sigma}_{i}$ is equivalent to that of the matrix pair $(A_{\rm\bf TT}(i), {\rm\bf col}\!\{ C_{\rm\bf T}(i),\; A_{\rm\bf ST}(i)\})$. Hence, the results of {Lemma 5} imply that, in order to construct an observable networked system, each of its subsystems should be observable.

To illustrate these theoretical results, we adopt a system model used in \cite{zhou15}. Due to space considerations, its parameters are not included. Interested readers are recommended to refer to \cite{zhou15} for details.

\hspace*{-0.45cm}{\bf Example I.} Consider the system of the first numerical example of \cite{zhou15}. Observability is checked for each subsystem respectively with external outputs {\it only} and with {\it both} the external and the internal outputs. The singular values of the associated observability matrices are given in Table I.

From these values, it is clear that except Subsystem ${\rm\bf\Sigma}_{1}$, {the other two subsystems are always observable, no matter they use  only their external outputs, that is, $y(t,i)$, or use both their internal and external outputs, that is, $z(t,i)$ and $y(t,i)$.} In addition, when both external and internal outputs are available, the first subsystem is also observable. The results of \cite{zhou15}, however, show that when the matrix $\Phi$ is utilized as the SCM, the overall system is unobservable. This confirms that observability of each subsystem can not guarantee that the overall system also has this property. On the other hand, when the SCM is modified to the matrix $\bar{\Phi}$, \cite{zhou15} shows that the overall system becomes observable. As $8.4218\times 10^{-10}$ is very close to zero and significantly less than one-thousandth of $1.1603$, it is reasonable to declare that the first subsystem is unobservable, at least very close to unobservable, provided that only its external outputs are available \cite{pzb14,zdg,zhou15}. Hence, appropriate subsystem connections can make the states of a subsystem observable that is unobservable with only its own external outputs.   \hspace{\fill}$\Diamond$

Note that observability of the matrix pair $(A_{\rm\bf TT}(i), {\rm\bf col}\!\{ C_{\rm\bf T}(i),\; A_{\rm\bf ST}(i)\})$ is not equivalent to that of the matrix pair $(A_{\rm\bf TT}(i), \; C_{\rm\bf T}(i))$. In fact, from Lemma 1, it is clear that if the matrix pair $(A_{\rm\bf TT}(i),\; C_{\rm\bf T}(i) )$ is observable, then {the matrix pair $(A_{\rm\bf TT}(i), {\rm\bf col}\!\{ C_{\rm\bf T}(i),\; A_{\rm\bf ST}(i)\})$ is also observable}; but the converse is in general not true. Results of Lemma 5 therefore also imply that even when there exist subsystems that are not observable through {\it only} their own external outputs, the whole networked system may still be observable by means of subsystem connections.

It is worthwhile to note that while similar results have been observed in \cite{zz17} for system controllability, the conclusions there depend on the SCM $\Phi$. This makes them difficult to be applied in constructing a controllable networked system, as an appropriate subsystem connection is usually not known before system designs. On the other hand, note that ${\rm\bf col}\!\{\lambda I_{M_{\rm\bf x}}\!\!-\!\!A_{\rm\bf TT},\; -C_{\rm\bf T},\; -\Phi A_{\rm\bf ST}\}
\!\!=\!\!{\rm\bf diag}\{I_{M_{\rm\bf x}},$ $-I_{M_{\rm\bf y}},\; -\Phi\}{\rm\bf col}\!\{\lambda I_{M_{\rm\bf x}}\!\!-\!\!A_{\rm\bf TT},\; C_{\rm\bf T},\; A_{\rm\bf ST}\}$. This means that in order to guarantee that the matrix ${\rm\bf col}\!\{\lambda I_{M_{\rm\bf x}}\!\!-\!\!A_{\rm\bf TT},\; -C_{\rm\bf T},\; -\Phi A_{\rm\bf ST}\}$ is of FCR, it is necessary that the matrix ${\rm\bf col}\!\{\lambda I_{M_{\rm\bf x}}\!\!-\!\!A_{\rm\bf TT},\; C_{\rm\bf T},\; A_{\rm\bf ST}\}$ is. Based on these observations and similar arguments as those in the proof of {Lemma 5}, it can be shown that the conclusions of {Lemma 5} are in fact valid for an arbitrary SCM $\Phi$.

To establish a relation between system observability and its subsystem out-degrees, define transfer function matrices $G^{[1]}(\lambda)$ and $G^{[2]}(\lambda)$ respectively as
$G^{[1]}(\lambda)\!=\!{\rm\bf diag}\{G^{[1]}_{i}(\lambda)|_{i=1}^{N}\}$ and $G^{[2]}(\lambda)\!=\!{\rm\bf diag}\{G^{[2]}_{i}(\lambda)|_{i=1}^{N}\}$, in which
$G^{[1]}_{i}(\lambda)\!=\!C_{\rm\bf S}(i)+C_{\rm\bf T}(i)[\lambda I_{m_{{\rm\bf x}i}}\!-\!A_{\rm\bf TT}(i)]^{\!-1}A_{\rm\bf TS}(i)$,
$G^{[2]}_{i}(\lambda)\!=\!A_{\rm\bf SS}(i)\!+\!A_{\rm\bf ST}(i)[\lambda I_{m_{{\rm\bf x}i}}\!-\!A_{\rm\bf TT}(i)]^{\!-1}A_{\rm\bf TS}(i)$ for each $i=1,2,\cdots,N$.
From the block diagonal structure of $G^{[1]}(\lambda)$, it is clear that it is of FCNR if and only if each of $G^{[1]}_{i}(\lambda)$, $i\in\{1,2,\cdots,N\}$,  is.

Assume that $G^{[1]}(\lambda)$ and $G^{[1]}_{i}(\lambda)$ have respectively $m$ and $m_{i}$ distinctive transmission zeros. Then, under the condition that $G^{[1]}(\lambda)$ is of FCNR, it is obvious from Lemma 2 and
$G^{[1]}(\lambda)\!=\!{\rm\bf diag}\{G^{[1]}_{i}(\lambda)|_{i=1}^{N}\}$ that, for each $i=1,\cdots,N$, every transmission zero of $G^{[1]}_{i}(\lambda)$ is also a transmission zero of $G^{[1]}(\lambda)$. As argued in \cite{zhou15}, we generally only have that $\max_{1\leq i\leq N}m_{i}\!\leq\! m\!\leq\! \sum_{i=1}^{N}m_{i}$. Moreover, for each of the transmission zeros of $G^{[1]}(\lambda)$, there exists at least one integer $i$ belonging to the set $\{\:1,\; 2,\;\cdots,\; N\:\}$, such that it is also a transmission zero of $G_{i}^{[1]}(\lambda)$.

Let $\lambda_{0}^{[k]}$ denote the $k$-th transmission zero of $G^{[1]}(\lambda)$, $k=1,2,\cdots,m$. Assume that in the set
$\{G_{1}^{[1]}(\lambda),\;G_{2}^{[1]}(\lambda),\;\cdots,$ $G_{N}^{[1]}(\lambda)\}$, there are $s^{[k]}$ transfer function matrices which have this transmission zero. Denote them by $G^{[1]}_{k(s)}(\lambda)$, $s=1,\cdots,s^{[k]}$. Clearly, both $s^{[k]}$ and $k(s)$ belong to the set $\{1,\;2,\;\cdots,\;N\}$. As in \cite{zhou15}, it is assumed, without any loss of generality, that $k(1)<k(2)<\cdots<k(s^{[k]})$. Let ${Y}^{[k]}_{s}$ denote the matrix constructed from  a set of linear independent vectors that span the null space of $G^{[1]}_{k(s)}(\lambda_{0}^{[k]})$, and $p(k,s)$ the dimension of this null space. Obviously, the matrix ${Y}^{[k]}_{s}$ is of FCR, which further leads to that the matrix ${Y}^{[k]H}_{s}{Y}^{[k]}_{s}$ is positive definite. Hence, the matrix $\Gamma_{s}^{[k]}$ is well defined for each $s=1,2,\cdots,s^{[k]}$ and each $k=1,2,\cdots,m$, which has the following definition
\begin{equation}
\Gamma_{s}^{[k]}=G_{k(s)}^{[2]}(\lambda_{0}^{[k]}){Y}^{[k]}_{s}\left({Y}^{[k]H}_{s}{Y}^{[k]}_{s}\right)^{-1/2}
\label{eqn:17}
\end{equation}

Using these matrices, the following conclusions are derived, which give a sufficient condition for the observability of the networked system $\rm\bf\Sigma$. Their proof is deferred to the appendix.

\hspace*{-0.45cm}{\bf Theorem 1.} Assume that all $G_{i}^{[1]}(\lambda)|_{i=1}^{N}$ are of FCNR. Let $\{\lambda_{0}^{[k]}|_{k=1}^{m}\}$ denote the set of distinctive transmission zero of $G^{[1]}(\lambda)$. If the matrix $\Theta$ satisfies simultaneously the following inequality
\begin{equation}
I_{p(k,s)}-\Gamma_{s}^{[k]H}\Theta^{2}(k(s))\Gamma_{s}^{[k]} >0
\label{eqn:24}
\end{equation}
or
\begin{equation}
I_{p(k,s)}-\Gamma_{s}^{[k]H}\Theta^{2}(k(s))\Gamma_{s}^{[k]} <0
\label{eqn:26}
\end{equation}
for each $s=1,2,\cdots, s^{[k]}$ and $k=1,2,\cdots,m$, then the dynamic system $\rm\bf\Sigma$ is observable.    \hspace{\fill}$\Diamond$

Compared with the results reported in \cite{zhou15}, the conditions of {Theorem 1} are only sufficient. On the other hand, these conditions can be verified individually for each subsystem and therefore have a much lower computational complexity, and the computation results are generally more numerically reliable. {In particular, in the above conditions, the dimension of the involved matrix is $p(k,s)\times p(k,s)$, while that in \cite{zhou15} is
$\sum_{i=1}^{N}m_{{\rm\bf v}i}\times \sum_{i=1}^{s^{[k]}}p(k,i)$. Obviously, the latter is usually significantly greater than the former for a large scale system, which is less attractive from the viewpoint of computations.}

Note that the matrix $\Theta$ is closely related to the out-degrees of the networked system $\rm\bf\Sigma$. {Theorem 1} in fact establishes a relation between the observability of a networked system and its subsystem out-degrees. {This theorem, together with the following Theorem 2, which is the counterpart of this theorem in controllability verifications, are essential in solving the minimal local input/output problem described in Section II. The details are given in the next section.}

\hspace*{-0.45cm}{\bf Remark I.} Note that for each $j=1,2,\cdots,N$, $\Theta(j)\geq I_{m_{{\rm\bf z}j}}$ from its definition. It can be easily understood that if there is an integer pair $(k,s)$ with $k\in\{1,2,\cdots,m\}$ and $s\in\{1,2,\cdots,s^{[k]}\}$, such that the associated matrix $\Gamma_{s}^{[k]}$ is not of FCR, then for all the SCM $\Phi$, the associated inequality $I_{p(k,s)}-\Gamma_{s}^{[k]H}\Theta^{2}(k(s))\Gamma_{s}^{[k]}<0$ can not be satisfied. Hence, to satisfy the conditions of {Theorem 1}, one possible approach is to meet the inequality $I_{p(k,s)}-\Gamma_{s}^{[k]H}\Theta^{2}(k(s))\Gamma_{s}^{[k]}>0$. This might be achieved by reducing the number of subsystems that an internal output straightforwardly affects. These observations mean that under such a situation, sparse subsystem connections might be helpful to make a networked system observable.

On the contrary, if for each $s=1,2,\cdots, s^{[k]}$ and each $k=1,2,\cdots,m$, the associated matrix $\Gamma_{s}^{[k]}$ is always of FCR, then the minimal eigenvalue of the matrix $\Gamma_{s}^{[k]H}\Theta^{2}(k(s))\Gamma_{s}^{[k]}$ can be made large through increasing the number of subsystems that an internal output directly influences, which implies that the inequality $I_{p(k,s)}-\Gamma_{s}^{[k]H}\Theta^{2}(k(s))\Gamma_{s}^{[k]}<0$ might be satisfied through simply increasing the number of subsystem connections. That is, dense subsystem connections are appreciated from the viewpoint of system observability. \hspace{\fill}$\Diamond$

\hspace*{-0.45cm}{\bf Example II.} To illustrate the influences of subsystem out-degrees on the observability of a networked system, consider again the system of the first numerical example in \cite{zhou15}. Modify its SCM $\Phi$ into the following one,
\begin{eqnarray*}
\Phi\!\!&=&\!\!{\rm\bf col}\{[0\;  0\; 1\; 0\; 0\; 0],\;\;[0\; 0\; 0\; 0\; 1\; 0],\;\;
      [1\; 0\; 0\; 0\; 0\; 0],\\
& &\hspace*{0.6cm} [0\; 0\; 0\; 0\; 0\; 1],\;\;[0\; 0\; 0\; 0\; 0\; 0],\;\; [0\; 1\; 0\; 0\; 0\; 0]\}
\end{eqnarray*}
Note that the fifth row of the SCM has been replaced by a row with all elements being zero, which means that the out-degree of the subsystem ${\rm\bf\Sigma}_{2}$ is reduced from $2$ to $1$. With this SCM, the singular values of the observability matrix of the whole networked system become $6.2043\times 10^{-1}$, $8.3070\times 10^{-1}$, $1.3217$, $1.5944$, $3.3636$, $6.7586$, $1.1178\times 10^{1}$, $1.5611\times 10^{1}$ and $9.0003\times 10^{3}$. It can therefore be declared that the associated networked system is now observable.   \hspace{\fill}$\Diamond$

\hspace*{-0.45cm}{\bf Remark II.} While the matrix $Y_{s}^{[k]}$ is not unique for each integer pair $(k,\;s)$, its selection does not have any influences on the satisfaction of the conditions of Equations (\ref{eqn:24}) and (\ref{eqn:26}), {which can be straightforwardly proven from relations among different
basis vectors of a subspace.} \hspace{\fill}$\Diamond$

When controllability is to be investigated, by means of the duality between controllability and observability of a LTI system, which has already been adopted in \cite{zhou15}, similar results can be derived through completely the same arguments. More precisely, based on this duality and the state space model of the whole system given in \cite{zhou15}, it can be directly declared that when the networked system $\rm\bf\Sigma$ is well-posed, it is controllable if and only if for each complex scale $\lambda$, the following matrix valued polynomial $\bar{M}(\lambda)$ is of FRR \cite{zhou15,zz17}
\begin{displaymath}
\bar{M}(\lambda)=\left[\begin{array}{ccc}
\lambda I_{M_{\rm\bf x}}-A_{\rm\bf TT} & -B_{\rm\bf T} & -A_{\rm\bf TS}\Phi \\
-A_{\rm\bf ST} & -B_{\rm\bf S} & I_{M_{\rm\bf v}}- A_{\rm\bf SS}\Phi \end{array}\right]
\end{displaymath}
Note that the transpose of $\bar{M}(\lambda)$ has completely the same form as that of $M(\lambda)$. It is not out of imaginations that  necessary/sufficient conditions similar to those of Lemma 5 and Theorem 1 can be derived for controllability verifications of a networked system.

However, in order to achieve these conclusions, it appears necessary to assume that every column of the SCM $\Phi$ only has one nonzero element. While this condition can be satisfied in general through augmenting the subsystem internal output vectors $z(t,i)|_{i=1}^{N}$ with repeated elements, the augmentation usually violates an associated FCNR condition and therefore greatly restricts applicability of the associated results.

In this paper, we derive another necessary/sufficient condition for system controllability without that assumption.

For this purpose, define $\bar{G}^{[1]}(\lambda)$ and $\bar{G}^{[2]}(\lambda)$ respectively as
$\bar{G}^{[1]}(\lambda)={\rm\bf diag}\{\bar{G}^{[1]}_{i}(\lambda)|_{i=1}^{N}\}$ and
$\bar{G}^{[2]}(\lambda)={\rm\bf diag}\{$ $\bar{G}^{[2]}_{i}(\lambda)|_{i=1}^{N}\}$, in which
$\bar{G}^{[1]}_{i}(\lambda)=B^{T}_{\rm\bf S}(i)+B^{T}_{\rm\bf T}(i)[\lambda I_{m_{\rm\bf Ti}}-$ $A^{T}_{\rm\bf TT}(i)]^{-1}A^{T}_{\rm\bf ST}(i)$ and
$\bar{G}^{[2]}_{i}(\lambda)=(G^{[2]}_{i}(\lambda))^{T}$.  Assume that $\bar{G}^{[1]}(\lambda)$ has $\bar{m}$ distinctive transmission zeros which are denoted by
$\bar{\lambda}_{0}^{[k]}|_{k=1}^{\bar{m}}$. Moreover, let $\bar{G}^{[1]}_{\bar{k}(s)}(\lambda)|_{s=1}^{\bar{s}^{[k]}}$ represent the transfer function matrices that have  $\bar{\lambda}_{0}^{[k]}$ as its transmission zero, and $\bar{k}(1)<\bar{k}(2)<\cdots<\bar{k}(\bar{s}^{[k]})$. Furthermore, let $\bar{p}(k,s)$ denote the dimension of
the null space of the matrix $\bar{G}^{[1]}_{\bar{k}(s)}(\bar{\lambda}_{0}^{[k]})$, and $\bar{Y}^{[k]}_{s}$ the matrix constructed from a set of linear independent vectors that span this null space. Define a matrix
$\bar{\Gamma}_{s}^{[k]}$ as
\begin{equation}
\bar{\Gamma}_{s}^{[k]}=\bar{G}_{\bar{k}(s)}^{[2]}(\bar{\lambda}_{0}^{[k]})\bar{Y}^{[k]}_{s}\left(\bar{Y}^{[k]H}_{s}\Theta^{-2}(\bar{k}(s))\bar{Y}^{[k]}_{s}\right)^{-1/2}
\label{eqn:30}
\end{equation}
Then, we have the following results, whose proof is included in the appendix.

\hspace*{-0.45cm}{\bf Theorem 2.} Assume that $\bar{G}^{[1]}(\lambda)$ is of FCNR. Then, System $\rm\bf\Sigma$ is controllable, only when the matrix pair $(A_{\rm\bf TT}(i),\;[B_{\rm\bf T}(i)\;\; A_{\rm\bf TS}(i)])$ is controllable for every $i=1,2,\cdots,N$. Moreover, if for each integer pair $(k,s)$ with $k\in\{1,2,\cdots,\bar{m}\}$ and $s\in\{1,2,\cdots,\bar{s}^{[k]}\}$, the following matrix inequality is satisfied,
\begin{equation}
I_{\bar{p}(k,s)}-\bar{\Gamma}_{s}^{[k]H}\bar{\Gamma}_{s}^{[k]}>0
\label{eqn:42}
\end{equation}
{then this system is controllable.}    \hspace{\fill}$\Diamond$

It is interesting to notice that while the necessary condition of Theorem 2 is  dual to that of Lemma 5, its sufficient condition differs significantly from that of Theorem 1. Moreover, their proofs are also not completely dual to each other. These are due to that in order to apply the duality between controllability and observability, it is necessary that the SCM $\Phi$ satisfies the condition that $\Phi\Phi^{T}$ is a diagonal matrix, which can not be met in general.

{From the definition of the matrix $\bar{\Gamma}_{s}^{[k]}$, careful comparisons between Equations (\ref{eqn:42}) and (\ref{eqn:24}) show that, some qualitative relations exist between in-degrees and controllability of a networked system, which are similar to those between its out-degrees and observability.}

\section{Minimal Local Input/Output Selection for a Networked System}

For a networked system, it is often interesting to know how many sensors are required to monitor its states, as well as how many actuators are required to maneuver its states \cite{pka16,pzb14,scl15,siljak78,trpj16,wj01}. Recall that in order to reconstruct the states of a system from measured input-output data, it is necessary that the system is observable. Moreover, controllability is necessary for a system to perform satisfactorily \cite{zdg,zhou15}. In this section, we investigate the minimal number of outputs/inputs required for {each subsystem} to guarantee the observability/controllability of the {whole} networked system, that is, the problem described in Section II, using the results of Section III.

\vspace{-0.0cm}
\begin{table*}[!ht]
\vspace*{-0.00cm}
\caption{\vspace*{-0.00cm} Singular Values of the Observability Matrices of the Overall System}
\begin{center}
\vspace{-0.2cm}
\begin{tabular}{|p{2.0cm}|p{1.8cm}|p{1.8cm}|p{1.8cm}|p{1.8cm}|p{1.8cm}|}
\hline
$k=1$ & $\kappa=0.98$ & $\kappa=0.96$  & $\kappa=0.94$ & $\kappa=0.92$ & $\kappa=0.90$ \\
\hline
\hline
  $1.9362\times 10^{-15}$ & $1.8098\times 10^{-1}$ &  $4.2275\times 10^{-1}$ & $4.9783\times 10^{-1}$ & $5.1380\times 10^{-1}$ & $5.2634\times 10^{-1}$   \\ \hline
  $4.4421\times 10^{-1}$ & $4.6358\times 10^{-1}$ &  $4.8757\times 10^{-1}$ & $7.6579\times 10^{-1}$ & $1.0373$ & $1.1071$   \\ \hline
  $1.0277$ & $1.0455$ &  $1.0744$ & $1.1353$ & $1.4116$ & $1.4448$   \\ \hline
  $1.4593$ & $1.4555$ &  $1.4524$ & $1.4509$ & $1.4803$ & $1.7099$   \\ \hline
  $2.1211$ & $2.0690$ &  $2.0268$ & $2.0120$ & $2.1752$ & $3.6578$   \\ \hline
  $4.0509$ & $4.0360$ &  $4.0194$ & $4.0057$ & $4.0171$ & $6.1553$   \\ \hline
  $7.3917$ & $7.4216$ &  $7.4645$ & $7.5180$ & $7.5792$ & $7.6485$   \\ \hline
  $1.1216\times 10^{1}$ & $1.1198\times 10^{1}$ &  $1.1187\times 10^{1}$ & $1.1183\times 10^{1}$ & $1.1187\times 10^{1}$ & $1.1224\times 10^{1}$   \\ \hline
  $2.1032\times 10^{1}$ & $2.0932\times 10^{1}$ &  $2.0822\times 10^{1}$ & $2.0703\times 10^{1}$ & $2.0577\times 10^{1}$ & $2.0454\times 10^{1}$   \\ \hline
\end{tabular}
\vspace{-0.5cm}
\end{center}
\end{table*}

The following theorem gives an answer to this minimal input/output problem. Its proof is provided in the appendix.

\hspace*{-0.45cm}{\bf Theorem 3.} Let $p_{max}(i)$ denote the maximum geometric multiplicity of the matrix $A_{\rm\bf TT}(i)$, $i=1,2,\cdots,N$. Then, an observable networked system $\rm\bf\Sigma$ can be constructed with local external outputs, if and only if
\begin{displaymath}
m_{{\rm\bf y}i}+m_{{\rm\bf z}i}\geq p_{max}(i),\hspace{0.5cm} \forall i\in\{\:1,\;2,\;\cdots,\;N\:\}
\end{displaymath}
Moreover, a controllable networked system $\rm\bf\Sigma$ can be constructed with local actuators, if and only if
\begin{displaymath}
m_{{\rm\bf u}i}+m_{{\rm\bf v}i}\geq p_{max}(i),\hspace{0.5cm} \forall i\in\{\:1,\;2,\;\cdots,\;N\:\}
\end{displaymath}
{Here, $m_{{\rm\bf *}i}$ stands for the dimension of the column vector $*(t,i)$ with $*=u,\;v,\;y,\;z$.  \hspace{\fill}$\Diamond$}

\hspace*{-0.45cm}{\bf Remark III.} This theorem reveals that in order to reduce the required number of external inputs/outputs, it is better to design a subsystem with its STM having distinctive eigenvalues. This is in a good agreement with the results on a lumped system reported in \cite{sm68,yzdwl13,zhou17}. \hspace{\fill}$\Diamond$

\hspace*{-0.45cm}{\bf Corollary 1.} In order to be able to build a controllable/observable networked system from several subsystems, it is necessary and sufficient that each subsystem is controllable/observable.

\hspace*{-0.45cm}{\bf Proof:} This is an immediate result of Lemmas 4 and 5, together with Theorem 3. \hspace{\fill}$\Diamond$

Note that the matrices $A_{\rm\bf ST}(i)$, $A_{\rm\bf SS}(i)$ and $A_{\rm\bf TS}(i)$ represent connection strengthes among subsystems of the system $\rm\bf\Sigma$. The bigger the magnitude of the elements of these matrices is, the tighter the subsystems are connected \cite{zhou15}. On the other hand, it is clear from the proof of Theorem 3 that when each subsystem is observable/controllable, through reducing subsystem connection strengthes, it is always possible to construct an observable/controllable networked system. In the extreme situation, when all the subsystems are disconnected, the networked system becomes a collection of isolated individual observable/controllable subsystems, which is obviously observable/controllable.

On the other hand, when these matrices are appropriately selected such that the corresponding matrices $\Gamma_{s}^{[k]}$ is of FCR for each integer pair $(k,\;s)$, it can be easily seen from Equation (\ref{eqn:26}) that through increasing magnitudes of the elements of these matrices, that is, through increasing subsystem connection strengthes, it is also possible to build an observable networked system using observable subsystems. Similar conclusions can be obtained for building a controllable networked system, by means of the duality between observability and controllability.

Combing together the results of Lemma 5 and Theorems 1 and 3, a prototypical algorithm can be constructed for building an observable networked system.

\hspace*{-0.45cm}{\rm\bf Algorithm for Constructing an Observable System:}
\begin{itemize}
\item Compute the maximum geometric multiplicity $p_{max}(i)$ for each STM $A_{\rm\bf TT}(i)$, $i=1,2,\cdots,N$. Select a real matrix $C(i)$ for Subsystem ${\rm\bf\Sigma}_{i}$ which has at least $p_{max}(i)$ rows, such that the matrix pair $(A_{\rm\bf TT}(i),\;C(i))$ is observable.
\item Partition the matrix $C(i)$ as $C(i)={\rm\bf col}\!\{ C_{\rm\bf T}(i),\; A_{\rm\bf ST}(i)\}$. In this partition, it is preferable to make the number of the rows of the matrix $A_{\rm\bf ST}(i)$ as small as possible, in order to reduce communication costs among subsystems.
\item Construct an initial value for each of the subsystem matrices $\left.A_{\rm\bf SS}(i)\right|_{i=1}^{N}$. Select an appropriate factor $\kappa$ belonging to $(0,\;1)$.
\begin{enumerate}
  \item Verify whether or not the constructed networked system is observable. If the answer is positive, end the computations.
  \item If the answer is negative, replace each subsystem matrix $A_{\rm\bf SS}(i)$ respectively with $\kappa A_{\rm\bf SS}(i)$, $1\leq i\leq N$. Return to Step 1).
\end{enumerate}
\end{itemize}

In the above algorithm, it is also possible to select a positive $\kappa$ greater than $1$. This situation, however, can be included through selecting an initial $\left.A_{\rm\bf SS}(i)\right|_{i=1}^{N}$ with the absolute values of their elements being large. In addition, similar algorithms can be constructed for building a controllable networked system using controllable subsystems.

{Note that for a lumped system, a complete parametrization has been given for its output matrix $C$ in \cite{zhou17} that constructs an observable matrix pair $(A,\;C)$ with a prescribed STM $A$. It can be declared that in the above algorithm, the construction of a desirable matrix $C(i)$ is not a difficult task.}

{Note also that in the above algorithm, the factor $\kappa$ satisfies $0<\kappa<1$. It is obvious that the magnitude of each element in the matrix $A_{\rm\bf SS}(i)$, $i=1,2,\cdots,N$, monotonically decreases with the iterations. As each subsystem is observable from its constructions, it can be declared from the proof of Theorem 3 that, with the increment of iterations, an observable networked system will certainly be constructed. On the other hand, the computational complexity of the above algorithm depends mainly on observability verification of the constructed system, for which a method is suggested in \cite{zhou15} whose computational complexity increases linearly with the subsystem number.}

\hspace*{-0.45cm}{\bf Example III} To reveal influences of subsystem connection strengthes on system observability, the subsystem matrices $\left.A_{\rm\bf SS}(i)\right|_{i=1}^{3}$ of the first numerical example in \cite{zhou15} is multiplied by a factor $\kappa$. Table II shows the singular values of the observability matrix of the whole system with $\kappa$ respectively equal to $1$, $0.98$, $0.96$, $0.94$, $0.92$ and $0.90$\footnote[1]{The first two singular values in the first column of Table II are a little different from those reported in \cite{zhou15}. This difference might be caused by numerical computations.}. Obviously, with the decrement of the subsystem connection strengthes, the whole system becomes observable.

Note that the eigenvalues of the subsystem STM $\left.A_{\rm\bf TT}(i)\right|_{i=1}^{3}$ are respectively
$9.0391\times 10^{-1}\pm j2.6359\times 10^{-2}$, $1.0243$; $6.7874\times 10^{-1}$, $8.6742\times 10^{-1}$, $1.0588$ and $8.5354\times 10^{-1}$, $1.2027$, $1.3343$.
Theorem 3 implies that only one output is required for each subsystem to make the whole system observable. Such a system can really be constructed through removing some elements of the output vectors $z(t,i)$ and/or $y(t,i)$, $i=1,2,3$. The details are omitted due to space considerations.

However, when there are restrictions on the SCM $\Phi$, and/or on subsystem connection strengthes, which is often required in practical engineering  \cite{sbkkmpr11,siljak78,zhou15}, further efforts are still necessary to find the minimal number of inputs/outputs for each subsystem in the construction of a controllable/observable system.

\section{Concluding Remarks}

In this paper, we have discussed minimal local input/output selections for a networked system. Some relations are established among out/in-degrees, observability and controllability of a networked system. It is observed that to guarantee the observability/controllability of the whole system, each subsystem must be observable/controllable. Moreover, according to properties of subsystems, sparse or dense connections may be helpful in constructing a controllable/observable system. Furthermore, in order to be able to construct a controllable/observable networked system, it is necessary and sufficient that each subsystem is controllable/observable, and the number of inputs/outputs of every subsystem must not be smaller than the maximum geometric multiplicity of its state transition matrix.

\renewcommand{\theequation}{a.\arabic{equation}}
\setcounter{equation}{0}

\section*{Appendix: Proof of Some Technical Results}

\hspace*{-0.45cm}{\rm\bf Proof of Lemma 5:} Define matrix valued polynomials $M_{1}(\lambda)$ and $\hat{M}_{1}(\lambda)$ respectively as
\begin{equation}
M_{1}(\lambda)\!=\!\left[\!\!\begin{array}{c}
\lambda I_{M_{\rm\bf x}}\!-\!A_{\rm\bf TT}  \\
-C_{\rm\bf T}  \\
-\Phi A_{\rm\bf ST} \end{array}\!\!\right],\hspace{0.25cm}
\hat{M}_{1}(\lambda)\!=\!\left[\!\!\begin{array}{c}
\lambda I_{M_{\rm\bf x}}\!-\!A_{\rm\bf TT}  \\
C_{\rm\bf T}  \\
\Theta A_{\rm\bf ST} \end{array}\!\!\right]
\label{eqn:6}
\end{equation}
Assume that the system ${\rm\bf\Sigma}$ is observable. Then, according to Lemma 3, it is necessary that for every complex scalar $\lambda$, the $M(\lambda)$ of Equation (\ref{eqn:4}) is of FCR. From the definitions of $M(\lambda)$ and $M_{1}(\lambda)$, it is obvious that $M_{1}(\lambda)$ must be of FCR at every complex scale $\lambda$ also, which is equivalent to
\begin{equation}
M_{1}^{H}(\lambda)M_{1}(\lambda)>0    \label{eqn:7}
\end{equation}

On the basis of Equations (\ref{eqn:5}) and (\ref{eqn:6}), the following equality can be straightforwardly established for each $\lambda\in{\cal C}$,
\begin{eqnarray}
M_{1}^{H}(\lambda)M_{1}(\lambda)\!\!&=&\!\!\left(\lambda I_{M_{\rm\bf x}}-A_{\rm\bf TT}\right)^{H}\left(\lambda I_{M_{\rm\bf x}}-A_{\rm\bf TT}\right)
+ \nonumber\\
& & \hspace*{2cm} C_{\rm\bf T}^{T}C_{\rm\bf T}+A_{\rm\bf ST}^{T}\Phi^{T} \Phi A_{\rm\bf ST}  \nonumber\\
&=&\!\! \left(\lambda I_{M_{\rm\bf x}}-A_{\rm\bf TT}\right)^{H}\left(\lambda I_{M_{\rm\bf x}}-A_{\rm\bf TT}\right)
+ \nonumber\\
& & \hspace*{2cm} C_{\rm\bf T}^{T}C_{\rm\bf T}+A_{\rm\bf ST}^{T}\Theta^{2} A_{\rm\bf ST}  \nonumber\\
&=&\!\! \hat{M}_{1}^{H}(\lambda)\hat{M}_{1}(\lambda)
\label{eqn:8}
\end{eqnarray}
It can therefore be declared that, in order to guarantee the observability of the system ${\rm\bf\Sigma}$, it is necessary that $\hat{M}_{1}(\lambda)$ is of FCR at each complex scale $\lambda$.

From the block diagonal structure of the matrices $A_{\rm\bf TT}$, $A_{\rm\bf ST}$ and $C_{\rm\bf T}$, as well as Equation (\ref{eqn:5}), it is obvious that
\begin{equation}
\hat{M}_{1}(\lambda)=\left[\begin{array}{c}
{\rm\bf diag}\!\left\{\! \lambda I_{m_{{\rm\bf x}i}}-A_{\rm\bf TT}(i)|_{i=1}^{N}\right\}  \\
{\rm\bf diag}\!\left\{\! C_{\rm\bf T}(i)|_{i=1}^{N}\right\}  \\
{\rm\bf diag}\!\left\{\! \Theta(i) A_{\rm\bf ST}(i)|_{i=1}^{N}\right\} \end{array}\right]     \label{eqn:9}
\end{equation}
Define matrix valued polynomials $\hat{M}_{1i}(\lambda)$ and $\tilde{M}_{1i}(\lambda)$ with $i=1,2,\cdots,N$ as
\begin{displaymath}
\hat{M}_{1i}(\lambda)\!=\!\!\left[\!\!\begin{array}{c}
\lambda I_{m_{{\rm\bf x}i}}\!-\!A_{\rm\bf TT}(i)  \\
C_{\rm\bf T}(i)  \\
\Theta(i) A_{\rm\bf ST}(i) \end{array}\!\!\right]\!,\hspace{0.25cm}
\tilde{M}_{1i}(\lambda)\!=\!\!\left[\!\!\begin{array}{c}
\lambda I_{m_{{\rm\bf x}i}}\!-\!A_{\rm\bf TT}(i)  \\
C_{\rm\bf T}(i)  \\
A_{\rm\bf ST}(i) \end{array}\!\!\right]
\end{displaymath}
Straightforward matrix manipulations show that for each fixed complex $\lambda$, the complex valued matrix $\hat{M}_{1}(\lambda)$ is of FCR, if and only if for each $i=1,2,\cdots,N$, the complex valued matrix $\hat{M}_{1i}(\lambda)$ is of FCR. Moreover, clearly from the definitions of $\hat{M}_{1i}(\lambda)$ and $\tilde{M}_{1i}(\lambda)$, we have that
\begin{equation}
\hat{M}_{1i}(\lambda)\!=\!{\rm\bf diag}\!\left\{\! I_{m_{{\rm\bf x}i}},\; I_{m_{{\rm\bf y}i}},\; \Theta(i)\right\}\tilde{M}_{1i}(\lambda)
\label{eqn:10}
\end{equation}
Note that the matrix $\Theta(i)$ is positive definite from its definition. It is clear that $\hat{M}_{1i}(\lambda)$ is of FCR at every complex scale $\lambda$, if and only if $\tilde{M}_{1i}(\lambda)$ is.

The proof can now be completed through a direct application of Lemma 1.  \hspace{\fill}$\Diamond$

\hspace*{-0.45cm}{\rm\bf Proof of Theorem 1:} From Lemma 3, it can be easily seen that System $\rm\bf\Sigma$ is observable, if and only if for each nonzero vector $x\in {\cal C}^{M_{\rm\bf x}+M_{\rm\bf v}}$, if there exists a $\lambda\in {\cal C}$, such that
\begin{equation}
\left[\begin{array}{cc}
\lambda I_{M_{\rm\bf x}}-A_{\rm\bf TT} & -A_{\rm\bf TS} \\
-C_{\rm\bf T} & -C_{\rm\bf S}  \end{array}\right]x=0      \label{eqn:11}
\end{equation}
then with the same complex number $\lambda$, the following inequality is valid
\begin{equation}
\left[-\Phi A_{\rm\bf ST} \;\; I_{M_{\rm\bf v}}-\Phi A_{\rm\bf SS} \right]x \neq 0      \label{eqn:12}
\end{equation}

Partition the vector $x$ as $x=\left[\: x_{1}^{T}\;\; x_{2}^{T}\:\right]$ in which $x_{1}\in {\cal C}^{M_{\rm\bf x}}$ and $x_{2}\in {\cal C}^{M_{\rm\bf v}}$. Then, according to Equation (\ref{eqn:11}), we have that
\begin{eqnarray}
& & \left[\lambda I_{M_{\rm\bf x}}-A_{\rm\bf TT}\right]x_{1}-A_{\rm\bf TS}x_{2}=0   \label{eqn:13} \\
& & C_{\rm\bf T}x_{1}+C_{\rm\bf S}x_{2}=0      \label{eqn:14}
\end{eqnarray}
When $\lambda$ is not an eigenvalue of the matrix $A_{\rm\bf TT}$, the matrix $\lambda I_{M_{\rm\bf x}}-A_{\rm\bf TT}$ is invertible. In this case, Equation (\ref{eqn:13}) implies that $x_{1}=\left[\lambda I_{M_{\rm\bf x}}-A_{\rm\bf TT}\right]^{-1}A_{\rm\bf TS}x_{2}$. Substitute this relation into Equations  (\ref{eqn:12}) and (\ref{eqn:14}), direct algebraic manipulations show that
\begin{eqnarray}
& & G^{[1]}(\lambda)x_{2}= 0 \label{eqn:15} \\
& & \left[I_{M_{\rm\bf v}}-\Phi G^{[2]}(\lambda)\right]x_{2} \neq 0      \label{eqn:16}
\end{eqnarray}
In these derivations, the definitions of the transfer function matrices $G^{[1]}(\lambda)$ and $G^{[2]}(\lambda)$ have been utilized.

When $\lambda$ is an eigenvalue of the matrix $A_{\rm\bf TT}$, a pseudo-inverse must be taken and the treatments are completely the same as those of \cite{zhou15,zz15}. In particular, note that the dimension of the matrix $A_{\rm\bf TT}$ is finite, which means that all its eigenvalues can only take an isolated value. Hence, there exists a $\varepsilon>0$ which in general depends on the value of $\lambda$, such that for each $\delta\in (-\varepsilon,\;\varepsilon)/\{0\}$, the matrix
$(\lambda-\delta) I_{M_{\rm\bf x}}-A_{\rm\bf TT}$ is invertible. These imply that the vector $x_{1}$ satisfying Equation (\ref{eqn:13}) can be formally expressed as
\begin{equation}
x_{1}=\lim_{\delta\rightarrow 0}\left[(\lambda-\delta) I_{M_{\rm\bf x}}-A_{\rm\bf TT}\right]^{-1}A_{\rm\bf TS}x_{2}
\end{equation}
Using this expression, conclusions can be obtained which are completely the same as those for the case when $\lambda$ is not an eigenvalue of the matrix $A_{\rm\bf TT}$.

Note that every $G_{i}^{[1]}(\lambda)$, $i=1, 2,\cdots, N$, is assumed to be of FCNR, and $G^{[1]}(\lambda)$ is block diagonal with its $i$-th diagonal block being $G^{[1]}_{i}(\lambda)$. It is obvious that $G^{[1]}(\lambda)$ is also of FCNR. It can therefore be declared from Lemma 2 and Equation (\ref{eqn:15}) that $\lambda$ is a transmission zero of $G^{[1]}(\lambda)$. These results imply that when $G_{i}^{[1]}(\lambda)|_{i=1}^{N}$ are of FCNR, verifications  of the conditions in Lemma 3 are necessary only for all the transmission zeros of $G^{[1]}(\lambda)$.

Assume that $\lambda=\lambda_{0}^{[k]}$. Then, according to the definition of the number $\lambda_{0}^{[k]}$, it is also a transmission zero of  $G_{k(s)}^{[1]}(\lambda)$, $s=1,2,\cdots, s^{[k]}$. Moreover, from the definition of the matrix $Y_{s}^{[k]}$, we have that for every nonzero complex valued vector $\alpha_{s}\in {\cal C}^{p(k,s)}$,
\begin{equation}
G_{k(s)}^{[1]}(\lambda_{0}^{[k]})Y_{s}^{[k]}\alpha_{s}=0      \label{eqn:37}
\end{equation}

Define a matrix $Y^{[k]}$ as
\begin{displaymath}
Y^{[k]}\!\!=\!\!\!\left[\!\!\begin{array}{ccc}
0_{M_{{\rm\bf v},k(1)\!-\!1}\!\times\! p(k,1)} & \cdots & 0_{M_{{\rm\bf v},k(s^{[k]})\!-\!1}\!\times\! p(k,s^{[k]})}    \\
Y_{1}^{[k]}  & \cdots & Y_{s^{[k]}}^{[k]} \\
0_{(M_{\rm\bf v}\!-\!M_{{\rm\bf v},k(1)})\!\times\! p(k,1)} & \cdots & 0_{(M_{\rm\bf v}\!-\!M_{{\rm\bf v},k(s^{[k]})})\!\times\! p(k,s^{[k]})}  \end{array}\!\!\!\right]
\end{displaymath}
Then, from the block diagonal structure of $G^{[1]}(\lambda)$ and Equation (\ref{eqn:37}), it can be directly proven that for each nonzero vector $x_{2}\in {\cal C}^{M_{\rm\bf v}}$ satisfying $G^{[1]}(\lambda_{0}^{[k]})x=0$, there exists one and only one nonzero $\alpha\in {\cal C}^{\sum_{j=1}^{s^{[k]}}p(k,j)}$, such that
\begin{equation}
x_{2}=Y^{[k]}\alpha    \label{eqn:18}
\end{equation}

On the other hand, based on the block diagonal structures of $G^{[2]}(\lambda)$ and the matrix $\Theta$, direct algebraic manipulations show that for each  complex valued vector $x_{2}$ satisfying Equation (\ref{eqn:18}), we have that
\begin{eqnarray}
& &\hspace*{-0.6cm}\Theta G^{[2]}(\lambda_{0}^{[k]})x_{2}\nonumber\\
& &\hspace*{-0.80cm}=\!{\rm\bf diag}\{\Theta(i)|_{i=1}^{N}\}{\rm\bf diag}\{G^{[2]}_{i}(\lambda)|_{i=1}^{N}\}Y^{[k]}\alpha \nonumber\\
& &\hspace*{-0.80cm}=\!\!\left[\!\!\!\!\begin{array}{ccc}
0_{M_{{\rm\bf v},k(1)-1}\times p(k,1)} & \cdots  & 0_{M_{{\rm\bf v},k(s^{[k]})-1}\times p(k,s^{[k]})}    \\
\Theta(k(1))G^{[2]}_{k(1)}\!(\lambda)\!Y_{1}^{[k]}  & \cdots  & \Theta(k(s^{[k]}))G^{[2]}_{k(s^{[k]})}\!(\lambda)\!Y_{s^{[k]}}^{[k]}  \\
0_{(M_{\rm\bf v}\!-\!M_{{\rm\bf v},k(1)})\!\times\! p(k,1)} & \cdots  & 0_{(M_{\rm\bf v}\!-\!M_{{\rm\bf v},k(s^{[k]})})\!\times\! p(k,s^{[k]})}
\!\!\!\!\!\!\!\!\end{array}\right]\!\!\alpha  \nonumber\\
\label{eqn:19}
\end{eqnarray}

\vspace{-0.4cm}Hence,
\begin{equation}
x_{2}^{H}x_{2}=\alpha^{H}{\rm\bf diag}\left\{ Y_{j}^{[k]H}Y_{j}^{[k]}|_{j=1}^{s^{[k]}} \right\}\alpha    \label{eqn:20}
\end{equation}
Moreover, from Equation (\ref{eqn:5}), we have that
\begin{eqnarray}
& &\!\!\!\!\left(\Phi G^{[2]}(\lambda_{0}^{[k]})x_{2}\right)^{H}\left(\Phi G^{[2]}(\lambda_{0}^{[k]})x_{2}\right)\nonumber\\
&=&\!\!\!\!x_{2}^{H} G^{[2]H}(\lambda_{0}^{[k]})\Theta^{2}G^{[2]}(\lambda_{0}^{[k]})x_{2} \nonumber \\
&=&\!\!\!\!\left(\Theta G^{[2]}(\lambda_{0}^{[k]})x_{2}\right)^{H}\left(\Theta G^{[2]}(\lambda_{0}^{[k]})x_{2}\right) \label{eqn:21}
\end{eqnarray}
Substitute the right hand side of Equation (\ref{eqn:19}) into that of Equation (\ref{eqn:21}), it can be directly proven that
\begin{eqnarray}
& &\!\!\!\!\left(\Phi G^{[2]}(\lambda_{0}^{[k]})x_{2}\right)^{H}\left(\Phi G^{[2]}(\lambda_{0}^{[k]})x_{2}\right)\nonumber\\
&=&\!\!\!\! \alpha^{H}{\rm\bf diag}\left\{ \left(\Theta(k(j)) G^{[2]}_{k(j)}(\lambda_{0}^{[k]}) Y_{j}^{[k]}\right)^{H}\times \right.\nonumber\\
& & \hspace*{0.8cm}\left.\left.\left(\Theta(k(j)) G^{[2]}_{k(j)}(\lambda_{0}^{[k]}) Y_{j}^{[k]}\right)\right|_{j=1}^{s^{[k]}} \right\}\alpha
\label{eqn:22}
\end{eqnarray}

Denote the vector ${\rm\bf diag}\{(Y_{j}^{[k]H}Y_{j}^{[k]})^{1/2}|_{j=1}^{s^{[k]}}\}\alpha$ by $\hat{\alpha}$. It can be declared from the FCR property of the matrices $Y_{j}^{[k]}|_{j=1}^{s^{[k]}}$ that the vector $\hat{\alpha}$ is not equal to zero if and only if the vector $\alpha$ is. On the other hand, from Equations (\ref{eqn:20}) and (\ref{eqn:22}), as well as the definitions of the matrices $\Gamma_{j}^{[k]}|_{j=1}^{s^{[k]}}$, straightforward algebraic manipulations show that
\begin{eqnarray}
\hspace*{-0.8cm}& &\!\!\!\!x_{2}^{H}x_{2}-\left(\Phi G^{[2]}(\lambda_{0}^{[k]})x_{2}\right)^{H}\left(\Phi G^{[2]}(\lambda_{0}^{[k]})x_{2}\right)\nonumber\\
\hspace*{-0.8cm} &=&\!\!\!\! \hat{\alpha}^{H}{\rm\bf diag}\left\{\left.I_{p(k,s)}-\Gamma_{s}^{[k]H}\Theta^{2}(k(s))\Gamma_{s}^{[k]} \right|_{s=1}^{s^{[k]}} \right\}\hat{\alpha}
\label{eqn:23}
\end{eqnarray}

Therefore, if the inequality of Equation (\ref{eqn:24}) is satisfied for each $s=1,2,\cdots,s^{[k]}$, then the matrix ${\rm\bf diag}\{I_{p(k,s)}-\Gamma_{s}^{[k]H}\Theta^{2}(k(s))\Gamma_{s}^{[k]}|_{s=1}^{s^{[k]}}\}$ is positive definite. This means that for an arbitrary nonzero complex vector $x_{2}$ satisfying Equation (\ref{eqn:15}), we have that
\begin{equation}
x_{2}^{H}x_{2}-\left(\Phi G^{[2]}(\lambda_{0}^{[k]})x_{2}\right)^{H}\left(\Phi G^{[2]}(\lambda_{0}^{[k]})x_{2}\right)> 0
\end{equation}
On the other hand, if for every $s\in\left\{1,2,\cdots,s^{[k]}\right\}$, the inequality of Equation (\ref{eqn:26}) is satisfied, then similar arguments show that for each nonzero complex vector $x_{2}$ satisfying Equation (\ref{eqn:15}), the following inequality is satisfied
\begin{equation}
x_{2}^{H}x_{2}-\left(\Phi G^{[2]}(\lambda_{0}^{[k]})x_{2}\right)^{H}\left(\Phi G^{[2]}(\lambda_{0}^{[k]})x_{2}\right)< 0
\end{equation}

Therefore, under both of these situations,
\begin{equation}
x_{2}\neq \Phi G^{[2]}(\lambda_{0}^{[k]})x_{2}
\label{eqn:25}
\end{equation}
Hence, $M(\lambda)$ is of FCR at each $\lambda=\lambda_{0}^{[k]}$. This completes the proof.  \hspace{\fill}$\Diamond$

\hspace*{-0.45cm}{\rm\bf Proof of Theorem 2:} To prove the condition for the necessity, assume that there exists a subsystem, denote it by ${\rm\bf\Sigma}_{i}$, such that the associated matrix pair $(A_{\rm\bf TT}(i),\;[B_{\rm\bf T}(i)\;\; A_{\rm\bf TS}(i)])$ is not controllable. Then, according to Lemma 1, there exist at least one $\lambda_{0}\in{\cal C}$ and one nonzero vector $x_{i}\in{\cal C}^{m_{{\rm\bf x}i}}$, such that
\begin{equation}
x_{i}^{H}\left[\lambda_{0}I_{m_{{\rm\bf x}i}}\!\!-\!\!A_{\rm\bf TT}(i)\;\;B_{\rm\bf T}(i)\;\; A_{\rm\bf TS}(i)\right]=0
\label{eqn:31}
\end{equation}
Define a $M_{\rm\bf x}$ dimensional vector $x$ as $x={\rm\bf col}\!\{0_{M_{{\rm\bf x},i-1}},\;x_{i},$ $\;0_{M_{{\rm\bf x}}-M_{{\rm\bf x},i}}\}$. Then, $x\neq 0$. Moreover, from Equation (\ref{eqn:31}) and the block diagonal structure of the matrices $A_{\rm\bf TT}$, $B_{\rm\bf T}$ and $A_{\rm\bf TS}$, direct matrix algebraic manipulations show that
\begin{equation}
x^{H}\left[\lambda_{0}I_{M_{{\rm\bf x}}}\!\!-\!\!A_{\rm\bf TT}\;\;-\!\!B_{\rm\bf T}\;\; -\!\!A_{\rm\bf TS}\right]=0
\label{eqn:32}
\end{equation}

Note that
\begin{eqnarray*}
& &\!\!\!\!\left[\lambda_{0}I_{M_{{\rm\bf x}}}\!\!-\!\!A_{\rm\bf TT} \;\;-\!\!B_{\rm\bf T}\;\; -\!\!A_{\rm\bf TS}\Phi\right] \\
&=&\!\!\!\!\left[\lambda_{0}I_{M_{{\rm\bf x}}}\!\!-\!\!A_{\rm\bf TT} \;\;-\!\!B_{\rm\bf T}\;\; -\!\!A_{\rm\bf TS}\right]{\rm\bf diag}\!\left\{\! I_{m_{{\rm\bf x}}},\;I_{m_{{\rm\bf u}}},\;\Phi\right\}
\end{eqnarray*}
We therefore have that the matrix $[\lambda_{0}I_{M_{{\rm\bf x}}}\!-\!A_{\rm\bf TT}\; -\!B_{\rm\bf T}\; -\!A_{\rm\bf TS}\Phi]$ can never be of FRR, no matter how the SCM $\Phi$ is designed. Hence, it can be claimed further from the definition of  $\bar{M}(\lambda)$ that it is also never of FRR at $\lambda=\lambda_{0}$. According to Lemma 1, System $\rm\bf\Sigma$ is not controllable.

To prove the condition for the sufficiency, note that $\bar{M}^{T}(\lambda)$ and ${M}(\lambda)$ have completely the same form. Similar arguments as those for the derivations of Equations (\ref{eqn:15}) and (\ref{eqn:16}) in the proof of {Theorem 1} show that, $\bar{M}(\lambda)$ is of FRR at each complex number $\lambda$, if and only if for each pair $(\lambda,\;x_{2})$ satisfying
\begin{equation}
\bar{G}^{[1]}(\lambda)x_{2}= 0 \label{eqn:33}
\end{equation}
in which $\lambda\in{\cal C}$, and $x_{2}\in{\cal C}^{{M_{\rm\bf z}}}$ and $x_{2}\neq 0$, the following inequality is satisfied
\begin{equation}
\left[I_{M_{\rm\bf z}}-\Phi^{T}\bar{G}^{[2]}(\lambda)\right]x_{2} \neq 0      \label{eqn:34}
\end{equation}

From the assumption that $\bar{G}^{[1]}(\lambda)$ is of FCNR and its block diagonal structure, as well as the definitions of the matrices $\bar{Y}_{s}^{[k]}|_{s=1}^{\bar{s}^{[k]}}$, it can be straightforwardly shown that every $\lambda$ satisfying Equation (\ref{eqn:33}) must be a transmission zero of $\bar{G}^{[1]}(\lambda)$. Moreover, all the nonzero $x_{2}$ satisfying Equation (\ref{eqn:33}) with $\lambda=\bar{\lambda}_{0}^{[k]}$ can be expressed as
\begin{equation}
x_{2}=\bar{Y}^{[k]}\alpha    \label{eqn:35}
\end{equation}
in which $\alpha$ is a nonzero $\sum_{s=1}^{\bar{s}^{[k]}}\bar{p}(k,s)$ dimensional complex vector and
\begin{displaymath}
\bar{Y}^{[k]}\!\!=\!\!\left[\!\!\begin{array}{ccc}
0_{M_{{\rm\bf z},\bar{k}(1)\!-\!1}\!\times\!\bar{p}(k,1)} & \cdots & 0_{M_{{\rm\bf z},\bar{k}(\bar{s}^{[k]})\!-\!1}\!\times\! \bar{p}(k,\bar{s}^{[k]})}    \\
\bar{Y}_{1}^{[k]}  & \cdots & \bar{Y}_{\bar{s}^{[k]}}^{[k]}  \\
0_{(M_{\rm\bf z}\!-\!M_{{\rm\bf z},\bar{k}(1)})\times \bar{p}(k,1)} & \cdots & 0_{(M_{\rm\bf z}\!-\!M_{{\rm\bf z},\bar{k}(\bar{s}^{[k]})})\!\times\! \bar{p}(k,\bar{s}^{[k]})} \end{array}\!\!\right]
\end{displaymath}

On the other hand, from Equation (\ref{eqn:5}) and singular value decompositions for a matrix \cite{hj91}, it can be declared that there exist a $U_{1}\in {\cal R}^{M_{\rm\bf v}\times M_{\rm\bf z}}$ and a $U_{2}\in {\cal R}^{M_{\rm\bf v}\times (M_{\rm\bf v}-M_{\rm\bf z})}$, such that
\begin{equation}
\Phi\!=\!U_{1}\Theta,\hspace{0.25cm} [U_{1}\; U_{2}]^{T}[U_{1}\; U_{2}]\!=\![U_{1}\; U_{2}][U_{1}\; U_{2}]^{T}\!=\!I_{M_{\rm\bf v}}
\label{eqn:36}
\end{equation}
Hence, for each $x_{2}$ satisfying Equation (\ref{eqn:35}), we have that
\begin{equation}
\left[\!I_{M_{\rm\bf z}}\!-\!\Phi^{T}\bar{G}^{[2]}(\lambda_{0}^{[k]})\!\right]x_{2}\!=\!\Theta
\left[\!\Theta^{-1}\bar{Y}^{[k]}\!-\!U_{1}^{T}\bar{G}^{[2]}(\lambda_{0}^{[k]})\bar{Y}^{[k]}\!\right]\alpha       \label{eqn:38}
\end{equation}
which means that $\left[\!I_{M_{\rm\bf z}}\!-\!\Phi^{T}\bar{G}^{[2]}(\lambda_{0}^{[k]})\!\right]x_{2}\neq 0$ if and only if
\begin{equation}
\left[\!\Theta^{-1}\bar{Y}^{[k]}\!-\!U_{1}^{T}\bar{G}^{[2]}(\lambda_{0}^{[k]})\bar{Y}^{[k]}\!\right]\alpha\neq 0        \label{eqn:39}
\end{equation}

Note that
\begin{equation}
\hspace*{-0.3cm}\left|\left|\Theta^{-1}\bar{Y}^{[k]}\alpha\right|\right|_{2}^{2}\!=\!\alpha^{H}\!
{\rm\bf diag}\!\left\{\!\left.\!\bar{Y}_{s}^{[k]H}\Theta^{-2}(\bar{k}(s))\bar{Y}_{s}^{[k]}\right|_{s=1}^{\bar{s}^{[k]}}\!\right\}\!\alpha
\label{eqn:40}
\end{equation}
Moreover, from Equation (\ref{eqn:36}), we have that $U_{1}U_{1}^{T}=I_{M_{\rm\bf v}}-U_{2}U_{2}^{T}\leq I_{M_{\rm\bf v}}$. Hence,
\begin{eqnarray}
& &\hspace*{-1.0cm}\left|\left|U_{1}^{T}\bar{G}^{[2]}(\lambda_{0}^{[k]})\bar{Y}^{[k]}\alpha\right|\right|_{2}^{2}\nonumber\\
& &\hspace*{-1.4cm}=\!\alpha^{H}\bar{Y}^{[k]H}\bar{G}^{[2]}(\lambda_{0}^{[k]H})U_{1}U_{1}^{T}\bar{G}^{[2]}(\lambda_{0}^{[k]})\bar{Y}^{[k]}\alpha \nonumber\\
& &\hspace*{-1.4cm}\leq\!\alpha^{H}\bar{Y}^{[k]H}\bar{G}^{[2]H}(\lambda_{0}^{[k]})\bar{G}^{[2]}(\lambda_{0}^{[k]})\bar{Y}^{[k]}\alpha \nonumber\\
& &\hspace*{-1.4cm}=\!\alpha^{H}\!{\rm\bf diag}\!\left\{\!\left.\bar{Y}_{s}^{[k]H}\bar{G}_{\bar{k}(s)}^{[2]H}(\lambda_{0}^{[k]})\bar{G}_{\bar{k}(s)}^{[2]}(\lambda_{0}^{[k]})\bar{Y}_{s}^{[k]}\right|_{s=1}^{\bar{s}^{[k]}} \!\right\}\!\alpha
\label{eqn:41}
\end{eqnarray}
which further leads to that
\begin{eqnarray}
& &\hspace*{-1.0cm}\left|\left|\Theta^{-1}\bar{Y}^{[k]}\alpha\right|\right|_{2}^{2}\!-\!\left|\left|U_{1}^{T}\bar{G}^{[2]}(\lambda_{0}^{[k]})\bar{Y}^{[k]}\alpha\right|\right|_{2}^{2}\nonumber\\
& &\hspace*{-1.4cm}\geq\!\alpha^{H}\!{\rm\bf diag}\!\left\{\!\left.\left(\bar{Y}_{s}^{[k]H}\Theta^{-2}(\bar{k}(s))\bar{Y}_{s}^{[k]}-\right.\right.\right.\nonumber\\
& & \hspace*{-0.0cm}\left.\left.\left.\bar{Y}_{s}^{[k]H}\bar{G}_{\bar{k}(s)}^{[2]H}(\lambda_{0}^{[k]})\bar{G}_{\bar{k}(s)}^{[2]}(\lambda_{0}^{[k]})\bar{Y}_{s}^{[k]}\right)\right|_{s=1}^{\bar{s}^{[k]}} \!\right\}\!\alpha   \nonumber\\
& &\hspace*{-1.4cm}=\!\hat{\alpha}^{H}\!{\rm\bf diag}\!\left\{\!\left.\left(I_{\bar{p}(k,s)}-\bar{\Gamma}_{s}^{[k]H}\bar{\Gamma}_{s}^{[k]}\right)\right|_{s=1}^{\bar{s}^{[k]}} \!\right\}\!\hat{\alpha}
\label{eqn:43}
\end{eqnarray}
in which $\hat{\alpha}={\rm\bf diag}\!\{(\bar{Y}_{s}^{[k]H}\Theta^{-2}(\bar{k}(s))\bar{Y}_{s}^{[k]})^{1/2}|_{s=1}^{\bar{s}^{[k]}} \!\}{\alpha}$.

Note that the matrix $\bar{Y}_{s}^{[k]H}\Theta^{-2}(\bar{k}(s))\bar{Y}_{s}^{[k]}$ is invertible for each feasible integer pair $(k,\;s)$. It is obvious that the vector $\alpha$ is nonzero if and only if the vector $\hat{\alpha}$ is. Therefore, if the condition of Equation (\ref{eqn:42}) is satisfied, then for any nonzero $\sum_{s=1}^{\bar{s}^{[k]}}\bar{p}(k,s)$ dimensional complex vector $\alpha$, we have that
\begin{equation}
\left|\left|\Theta^{-1}\bar{Y}^{[k]}\alpha\right|\right|_{2}^{2}\!-\!\left|\left|U_{1}^{T}\bar{G}^{[2]}(\lambda_{0}^{[k]})\bar{Y}^{[k]}\alpha\right|\right|_{2}^{2}>0
\end{equation}
Hence, the condition of Equation (\ref{eqn:39}) is satisfied, which means that the system $\rm\bf\Sigma$ is controllable. This completes the proof.  \hspace{\fill}$\Diamond$

\hspace*{-0.45cm}{\rm\bf Proof of Theorem 3:} From Theorem 1, we have that in order to guarantee the observability of the networked system $\rm\bf\Sigma$, it is necessary that for each $i=1,2,\cdots,N$, the matrix pair $(A_{\rm\bf TT}(i), [C^{T}_{\rm\bf T}(i)\; A^{T}_{\rm\bf ST}(i)]^{T})$ is observable. It can therefore be declared from Lemma 4 that to construct an observable $\rm\bf\Sigma$, it is necessary that $m_{{\rm\bf y}i}+m_{{\rm\bf z}i}\geq p_{max}(i)$.

Now, assume that $m_{{\rm\bf y}i}+m_{{\rm\bf z}i}= p_{max}(i)$ for every $1\leq i\leq N$. Then, according to Lemma 4, there always exists a matrix ${C}_{\rm\bf T}(i)$ and a matrix ${A}_{\rm\bf ST}(i)$ for each $i\in \{1,2,\cdots, N\}$, such that the matrix pair $(A_{\rm\bf TT}(i), [C^{T}_{\rm\bf T}(i)\; A^{T}_{\rm\bf ST}(i)]^{T})$ is observable.

Note that for an arbitrary real number $\kappa_{i}$, we have that
\begin{displaymath}
\left[\!\!\!\!\begin{array}{c} \lambda I_{m_{{\rm\bf x}i}}\!\!-\!\!A_{\rm\bf TT}(i) \\ {C}_{\rm\bf T}(i) \\ \kappa_{i}{A}_{\rm\bf ST}(i)\end{array}\!\!\!\!\right] \!\!\!=\!{\rm\bf diag}\!\!\left\{\!I_{m_{{\rm\bf x}i}}\!,\: I_{m_{{\rm\bf y}i}}\!,\:\kappa_{i}I_{m_{{\rm\bf z}i}}\!\right\}\!\!\!
\left[\!\!\!\!\begin{array}{c} \lambda I_{m_{{\rm\bf x}i}}\!\!-\!\!A_{\rm\bf TT}(i) \\ {C}_{\rm\bf T}(i) \\ {A}_{\rm\bf ST}(i) \end{array}\!\!\!\!\right]
\end{displaymath}
It is clear from Lemma 1 that observability of the matrix pair $(A_{\rm\bf TT}(i), {\rm\bf col}\!\{C_{\rm\bf T}(i),\; \kappa_{i}A_{\rm\bf ST}(i)\})$ is equivalent to that of the matrix pair $(A_{\rm\bf TT}(i), {\rm\bf col}\!\{C_{\rm\bf T}(i),\; A_{\rm\bf ST}(i)\})$, provided that $\kappa_{i}\neq 0$.

For each $j\in\{\:1,\;2,\;\cdots,\;N\:\}$, define a set ${\cal J}(j)$ as
\begin{displaymath}
{\cal J}(j)=\left\{\:(k,s)\:\left|\: k(s)=j,\; \begin{array}{l} s\in\{\:1,\;2,\;\cdots,\;s^{[k]}\:\}  \\
k\in \{\:1,\;2,\;\cdots,\;m\:\}\end{array} \right.\right\}
\end{displaymath}
That is, this set is associated with all the transmission zeros of $G^{[1]}(\lambda)$ that is also a transmission zero of $G_{j}^{[1]}(\lambda)$ with $j\in \{\:1,\;2,\;\cdots,\;N\:\}$. Then, obviously, the satisfaction of Equation (\ref{eqn:24}) can be equivalently expressed as
that for each $j=1,2,\cdots, N$, the following inequality
\begin{equation}
I_{p(k,s)}-\Gamma_{s}^{[k]H}\Theta^{2}(j)\Gamma_{s}^{[k]}>0
\label{eqn:27}
\end{equation}
is satisfied for every pair $(k,s)$ of the set ${\cal J}(j)$.

For a fixed SCM $\Phi$, define $\gamma_{i}$ as
\begin{equation}
\gamma_{i}\!=\!\max\!\!\left\{\!{\sigma}_{max}\!\!\left(\!\Theta(i){A}_{\rm\bf SS}(i)\!\right)\!,\hspace{0.05cm}\max_{(k,s)\in {\cal J}(i)}\!\!\!\!\!\!{\sigma}_{max}\!\!\left(\!\Theta(i)\Gamma_{s}^{[k]}\!\right)\!\!\right\}
\label{eqn:28}
\end{equation}
in which ${\sigma}_{max}(\cdot)$ stands for the maximal singular value of a matrix. Moreover, for each subsystem of System $\rm\bf\Sigma$, define matrices $\hat{A}_{\rm\bf ST}(i)$ and $\hat{A}_{\rm\bf SS}(i)$ respectively as
\begin{equation}
\hat{A}_{\rm\bf ST}(i)=\kappa_{i}{A}_{\rm\bf ST}(i), \hspace{0.5cm}
\hat{A}_{\rm\bf SS}(i)=\kappa_{i}{A}_{\rm\bf SS}(i)
\label{eqn:29}
\end{equation}
in which $\kappa_{i}$ is an arbitrary number belonging to $(0,\;1/\gamma_{i})$.

Using these two matrices, construct a new networked system $\hat{\rm\bf\Sigma}$ through simply replacing the system matrices $A_{\rm\bf ST}(i)$ and $A_{\rm\bf SS}(i)$ respectively by $\hat{A}_{\rm\bf ST}(i)$ and $\hat{A}_{\rm\bf SS}(i)$, while keeping the other system matrices unchanged. Moreover, define matrices $\hat{A}_{\rm\bf SS}$, $\hat{A}_{\rm\bf ST}$, etc., as well as transfer function matrices $\hat{G}^{[1]}(\lambda)$, $\hat{G}^{[2]}(\lambda)$, etc., respectively as their counterparts associated with System ${\rm\bf\Sigma}$.

Based on the block diagonal structure of the matrix $\hat{A}_{\rm\bf SS}$ and Equation (\ref{eqn:5}), it can be straightforwardly proven that $(\Phi \hat{A}_{\rm\bf SS})^{T}\!(\Phi \hat{A}_{\rm\bf SS})\!=\!{\rm\bf diag}\{\kappa_{i}^{2} A^{T}_{\rm\bf SS}(i)\Theta^{2}(i)A_{\rm\bf SS}(i)|_{i=1}^{N}\}$. Hence, it can be claimed from Equations (\ref{eqn:28}) and (\ref{eqn:29}) that
\begin{eqnarray}
{\sigma}_{max}\left(\Phi \hat{A}_{\rm\bf SS}\right)&=&\max_{1\leq i\leq N}\left\{
{\sigma}_{max}\left(\Theta(i) \hat{A}_{\rm\bf SS}(i)\right)\right\} \nonumber\\
&=&\max_{1\leq i\leq N}\left\{\kappa_{i}\times
{\sigma}_{max}\left(\Theta(i){A}_{\rm\bf SS}(i)\right)\right\} \nonumber\\
&<&  1
\end{eqnarray}
Note that the absolute value of each eigenvalue of a square matrix is not greater than its maximal singular value \cite{hj91}. It can therefore be declared that the matrix $I\!-\!\Phi\hat{A}_{\rm\bf SS}$ is invertible, and hence the re-constructed networked system $\hat{\rm\bf\Sigma}$ is well-posed.

On the other hand, note that in System $\hat{\rm\bf\Sigma}$, only the matrices $\hat{A}_{\rm\bf ST}(i)$ and $\hat{A}_{\rm\bf SS}(i)$ are different from those of System ${\rm\bf\Sigma}$. This implies that $G^{[1]}(\lambda)$ and $\hat{G}^{[1]}(\lambda)$, their transmission zeros, as well as the associated matrices $Y_{s}^{[k]}$, are completely the same. It can therefore be declared from the definition of the matrix $\Gamma_{s}^{[k]}$ that for each integer pair $(k,\;s)$ with $k\in\{\:1,\;2,\;\cdots,\;m\:\}$ and $s\in\{\:1,\;2,\;\cdots,\;s^{[k]}\:\}$, there certainly exists one and only one $j\in\{\:1,\;2,\;\cdots,\;N\:\}$, such that
the pair $(k,s)$ belongs to the set ${\cal J}(j)$. This further leads to that
\begin{equation}
\hat{\Gamma}_{s}^{[k]}=\kappa_{j}\Gamma_{s}^{[k]}
\end{equation}
Hence, we have from Equations (\ref{eqn:28}) and (\ref{eqn:29}) that
\begin{equation}
\sigma_{max}\left(\hat{\Gamma}_{s}^{[k]}\Theta(j)\right)=\kappa_{j}\sigma_{max}\left({\Gamma}_{s}^{[k]}\Theta(j)\right)< 1
\end{equation}
which further implies the satisfaction of the condition of Equation (\ref{eqn:27}) for each element of the set ${\cal J}(j)$ and each $j\in \{\:1,\;2,\;\cdots,\;N\:\}$, and hence the system $\hat{\rm\bf\Sigma}$ is observable.

The results on minimal input selection for system controllability can be established directly using duality between controllability and observability of a dynamic system, as well as the sufficient condition of Theorem 2.

This completes the proof. \hspace{\fill}$\Diamond$

\end{document}